\documentclass[12pt,a4paper, leqno]{article}
\usepackage{amsmath}
\usepackage{amssymb}
\usepackage{amsthm}
\usepackage{mathrsfs}
\usepackage[dvips]{graphicx}
\usepackage{psfrag}
\usepackage[hang,small,bf]{caption}
\usepackage{fullpage}
\newtheorem{theorem}{Theorem}[section]
\newtheorem{corollary}[theorem]{Corollary}

\newtheorem{lemma}[theorem]{Lemma}

\newtheorem{remark}[theorem]{Remark}
\numberwithin{equation}{section}
\newcommand{\G}{\mathcal{G}}
\newcommand{\E}{\mathcal{E}}
\newcommand{\C}{a}

\begin{document}
\title{Interlacement percolation on transient weighted graphs}
\author{A. Teixeira\footnote{Author's webpage: http://www.math.ethz.ch/\textasciitilde teixeira/}
 \\ {\small Departement Mathematik - ETH Z\"urich} \\ {\small R\"amistrasse, 101 HG G 47.2,} \\{\small CH-8092 Z\"urich, Switzerland} \\{\small\textbf{augusto.teixeira@math.ethz.ch}}}
\date{Submitted at June 24, 2008. Accepted at June 22, 2009.}
\maketitle

\begin{abstract}
In this article, we first extend the construction of random interlacements, introduced by A.S. Sznitman in \cite{sznitman}, to the more general setting of transient weighted graphs. We prove the Harris-FKG inequality for this model and analyze some of its properties on specific classes of graphs. For the case of non-amenable graphs, we prove that the critical value $u_*$ for the percolation of the vacant set is finite. We also prove that, once $\G$ satisfies the isoperimetric inequality $I \negmedspace S_6$ (see (\ref{eq:strongperim})), $u_*$ is positive for the product $\G \times \mathbb{Z}$ (where we endow $\mathbb{Z}$ with unit weights). When the graph under consideration is a tree, we are able to characterize the vacant cluster containing some fixed point in terms of a Bernoulli independent percolation process. For the specific case of regular trees, we obtain an explicit formula for the critical value $u_*$.
\end{abstract}

\section{Introduction}
\label{sec:intro}
The model of random interlacements was recently introduced by A.S. Sznitman in \cite{sznitman}. Its definition is motivated by the study of the trajectory performed by random walk on the discrete torus $(\mathbb{Z}/N\mathbb{Z})^d$, $d \geqslant 3$, or on the discrete cylinder $(\mathbb{Z} /N \mathbb{Z})^d \times \mathbb{Z}$ ,$d \geqslant 2$, when it runs up to times of respective order $N^d$ and $N^{2d}$, see \cite{benjamini}, \cite{sznitman_cil0}. In a heuristic sense, the random interlacements describe for these time scales the microscopic limiting ``texture in the bulk'' created by the walk, see \cite{david2}, \cite{sznitman_cil1}. In \cite{sznitman_cil3}, the random interlacements came as the main ingredient to improve the upper bound obtained in \cite{sznitman_cil0} for the asymptotic behavior of the disconnection time of the cylinder $(\mathbb{Z} /N \mathbb{Z})^d \times \mathbb{Z}$ by a random walk.

In this article, we first extend the model of random interlacements to the setting of transient weighted graphs, as suggested in \cite{sznitman}, Remark~1.4. Consider $\G = (V,\E)$ a graph composed of a countable set of vertices $V$ and a (non-oriented) set of edges $\E \subset \{ \mbox{\fontsize{10}{10} \selectfont $\{x,y\}$} \subset V;$ $x \neq y\}$. We assume $\G$ connected and provide it with a weight function $\C:V \times V \rightarrow \mathbb{R}_+$ (also called conductance), which is symmetric and such that $\C_{x,y} > 0$ if and only if $\{x,y\} \in \E$. When $\C_{x,y} = 1_{ \{x,y\} \in \E}$, we say that the graph is endowed with the canonical weights. A weight function $\C$ induces on $\G$ an irreducible, reversible Markov chain with transition probability given by $q(x,y) = \C_{x,y}/\mu_x$, where $\mu$ is the reversible measure for the chain, which is defined by $\mu_x = \sum_{\{x,y\} \in \E} \C_{x,y}$. Throughout this paper, we assume that the weighted graph under consideration induces a transient Markov chain.

Roughly speaking, the random interlacements are defined in terms of a Poisson point process on the space of doubly infinite trajectories in $V$ modulo time shift that visit points finitely often, see Section~\ref{sec:definition} for the precise definition. We will be interested in $\mathcal{I}^u$, the so-called interlacement at level $u \geqslant 0$, which is the union of the trace of the trajectories appearing in the above mentioned point process. The parameter $u$ controls the intensity of this process.

Although the precise definition of $\mathcal{I}^u$ is postponed to Section~\ref{sec:definition}, we now describe the law of the indicator function of the vacant set at level $u$, $\mathcal{V}^u = V \setminus \mathcal{I}^u$, regarded as a random element of $\{0,1\}^V$. As we show in Remark~\ref{prop:Q_u}, the law $Q^u$ that $\mathcal{V}^u$ induces on $(\{0,1\}^V, \mathcal{Y})$ is characterized by
\begin{equation}
\label{eq:Q_u}
Q^u[Y_x = 1 \text{ for all } x \in K] = \text{exp}(-u \cdot \text{cap}(K)), \text{ for all finite $K \subset V$},
\end{equation}
where $\text{cap}(K)$ stands for the capacity of $K$, see (\ref{eq:capacity}), and  $\mathcal{Y}$ is the $\sigma$-algebra generated by the canonical coordinate maps $(Y_x)_{x \in V}$ on $\{0,1\}^V$.

\vspace*{4mm}

We then prove in Theorem~\ref{th:harry}, a Harris-FKG type inequality for the law $Q^u$, answering a question of \cite{sznitman}, cf. Remark~1.6 1). More precisely, we show that for every pair of increasing random variables $f$ and $g$, see the beginning of Section~\ref{sec:definition} for the precise definition, with finite second moment with respect to $Q^u$, one has
\begin{equation}
\label{eq:harrisfkg}
\int fg \, dQ^u \geqslant \int f \, dQ^u \int g \, dQ^u.
\end{equation}

As a by-product, (\ref{eq:harrisfkg}) enables to define the critical value
\begin{equation}
\label{eq:critical}
u_* = \inf \{u \geqslant 0; \eta(x,u) = 0\} \in [0,\infty],
\end{equation}
for the percolation probability
\begin{equation}
\label{eq:eta}
\eta(u,x) = \mathbb{P}\left[\text{the cluster containing $x$ in $\mathcal{V}^u$ is infinite}\right],
\end{equation}
independently of the base point $x$. This is the content of Corollary~\ref{cor:uxu}. An important question is to determine whether $u_*$ is non-degenerate, i.e. $0<u_*<\infty$.

\vspace*{4mm}

In the remainder of the article, we derive some properties of random interlacements for specific classes of graphs.

We first derive two results concerning the non-degeneracy of the critical value $u_*$ under assumptions involving certain isoperimetric inequalities, $I \negmedspace S_d$, $d \geqslant 1$. Namely, according to \cite{woess} p. 40, a weighted graph $\G$ satisfies the isoperimetric inequality $I \negmedspace S_d$ if there exists $\kappa > 0$ such that
\begin{equation}
\label{eq:strongperim}
\mu(A)^{1-1/d} \leqslant \kappa \cdot \C(A), \text{ for every finite } A \in V,
\end{equation}
where
\begin{gather}
\label{eq:cA}
\mu(A) = \sum_{x \in A} \mu_x \text{ and}\\
\C(A) = \sum_{x \in A, y \in A^c} \C_{x,y}.
\end{gather}
When $d = \infty$ (with the convention $1-1/d = 1$), the graph is said to satisfy the strong isoperimetric inequality, or to be non-amenable.

Our first result concerning the non-degeneracy of $u^*$ states that, cf. Theorem~\ref{th:nonamenable},
\begin{equation}
\label{eq:statnonamenable}
\begin{split}
\text{$u_*$ is} & \text{ finite for non-amenable graphs with bounded}\\
& \text{degrees and weights bounded from below.}
\end{split}
\end{equation}

We then consider the critical value of product graphs of the type $\G \times \mathbb{Z}$. Random interlacements on such graphs are expected to be related to the investigation of the disconnection time of a discrete cylinder, see \cite{sznitman_how}.
In Theorem~\ref{th:isoper} we prove that
\begin{equation}
\label{eq:statis6}
\begin{split}
\text{the critical value $u_*$ of $\G \times \mathbb{Z}$ } & \text{is positive, if $\G$ satisfies the}\\
\text{isoperimetric inequality $I \negmedspace S_6$,} & \text{ has bounded degree and}\\
\text{weights bounded from a} & \text{bove and from below,}
\end{split}
\end{equation}
where we endow $\mathbb{Z}$ with the canonical weights and define the product weighted graph $\G \times \mathbb{Z}$ as in (\ref{eq:productweight}).

In particular, with these results, one concludes that when $\G$ as above is non-amenable, the critical value of $\G \times \mathbb{Z}$ is non-degenerate, see below (\ref{eq:eta}).


\vspace*{4mm}

We then consider the case where $\G$ is a transient tree of bounded degree. In general, the law $Q^u$ does not dominate nor is dominated by any non-degenerate Bernoulli i.i.d. variables, see Remark~1.1 of \cite{vladas}. However, when $\G$ is a tree of bounded degree, we show in Theorem~\ref{th:tree} that under the measure $Q^u$,
\begin{equation}
\label{eq:stattree}
\begin{split}
\text{the cluster containing a } & \text{fixed site $x \in V$ has the same law as}\\
\text{the cluster of $x$ under the } & \text{Bernoulli independent site percolation}\\
\text{process characterized } & \text{by $P[z \textnormal{ is open} \thinspace]$ $=$ $\textnormal{exp}(-u \cdot f_x(z))$,}
\end{split}
\end{equation}
where the function $f_x(z)$ is defined in (\ref{eq:bernoullicouple}). As a consequence, we conclude in Proposition~\ref{prop:Q_u} that
\begin{equation}
\label{eq:nondegtree}
\begin{split}
\text{w} & \text{hen $\G$ is a tree with degree, which is bounded,}\\
\text{an} & \text{d at least three, endowed with weights bounded}\\
\text{f} & \text{rom above and from below, then $0 < u_* < \infty$}.
\end{split}
\end{equation}

We also obtain an explicit formula for the critical value of regular trees of degree $d$:
\begin{equation}
\label{eq:criticaltreeintro}
u_* = \frac{d(d-1)\text{log}(d-1)}{(d-2)^2},
\end{equation}
see Corollary~\ref{cor:regular}. Interestingly this implies that $\mathbb{P}[x \in \mathcal{V}^{u_*}] = \frac{1}{d} (1 + o(1))$, as $d \rightarrow \infty$, cf Remark 5.3.

\vspace*{4mm}

We now give a rough description of the proofs of the main results in this article.

The construction of the intensity measure of the Poisson point process governing the random interlacements is the main step towards the definition of the process. It appears in Theorem~\ref{th:exist_nu}. Although we follow the argument of Theorem~1.1 of \cite{sznitman} for the case $V = \mathbb{Z}^d$, we present here the entire proof for the sake of completeness.

Concerning the Harris-FKG inequality, we cannot rely on the so-called FKG-Theorem (see \cite{liggett} Corollary 2.12 p. 78) to prove (\ref{eq:harrisfkg}). Indeed, the canonical condition (2.13) of \cite{liggett} p. 78 does not hold for the measure $Q^u$, see Remark~\ref{rem:harrisfinite} 2). Moreover, we also show in Remark~\ref{rem:harrisfinite} 1) that in general the conditional expectations of increasing functions on $\{0,1\}^V$, with respect to the $\sigma$-algebra generated by finitely many coordinates, are not necessarily increasing functions. This last feature prevents the use of the standard argument employed to prove the full Harris-FKG inequality once it holds for random variables depending only on the state of finitely many sites, see for instance \cite{grimmett} Theorem~2.4.

It is not clear how to obtain (\ref{eq:harrisfkg}) from the characterization of the law $Q^u$ given in (\ref{eq:Q_u}). Our proof relies instead on the construction of $\mathcal{I}^u$ in terms of the point process of interlacement trajectories, as in Section \ref{sec:definition}. Roughly speaking, we first transport the functions $f$ and $g$ to the space of point measures where this point process is defined. We then consider a sequence of $\sigma$-algebras $\mathcal{F}_n$ in this space, induced by an increasing sequence of finite sets $K_n$ exhausting $V$. Intuitively, these $\sigma$-algebras keep track of the behavior of interlacement trajectories of the point measure, which meet $K_n$, between their first and last visit to $K_n$. We then prove (\ref{eq:harrisfkg}) for $\mathcal{F}_n$-measurable functions using the FKG-Theorem, implying the desired result via a martingale convergence argument.


\vspace*{4mm}

For the proof of (\ref{eq:statnonamenable}), we rely on the known fact that for non-amenable graphs, the $L^2(\mu)$ norm of a compactly supported function $f$ on $V$ can be bounded in terms of the Dirichlet form of $f$
$$
\mathcal{D}(f,f) = \frac{1}{2} \sum_{x,y \in V} |f(x) - f(y)|^2 \C_{x,y},
$$
see \cite{woess} Theorem~10.3. This implies a linear lower bound for the capacity of a finite set in terms of its cardinality. This bound is then used to offset the growth of the number of self-avoiding paths of size $n$, in a Peierls-type argument leading to the finiteness of $u_*$.

\vspace*{4mm}

To prove (\ref{eq:statis6}) we use a renormalization argument that takes place on an isometric copy of the discrete upper half plane $\mathbb{Z}_+ \times \mathbb{Z}$, which can be found inside the graph $\G \times \mathbb{Z}$, cf. below (\ref{eq:taudistance}). We then employ bounds on the hitting probability of a point $x \in V$ for the random walk in $\G \times \mathbb{Z}$ in terms of the distance between $x$ and the starting point of the walk, which are consequences of classical results on isoperimetric inequalities, see for instance \cite{woess}, Theorem~14.3, p. 148. The proof of the positivity of $u_*$ then relies on a renormalization argument, which is adapted from \cite{sznitman}, Proposition~4.1.

\vspace*{4mm}

Theorem~\ref{th:tree} provides a characterization of the cluster of the vacant set containing a fixed site $x \in V$, in terms of a Bernoulli independent site percolation process, in the case where $\G$ is a tree, see (\ref{eq:stattree}). The rough strategy of the proof is to partition the space of doubly infinite trajectories modulo time shift, where the Poisson point process governing $\mathcal{I}^u$ is defined, into sets indexed by the vertices of the graph $V$. This partition induces on $V$ a Bernoulli independent site percolation process, and we can identify the corresponding component of $x$ with the connected component of $x$ in $\mathcal{V}^u$. We strongly use the fact that $\G$ is a tree to obtain this identification.

\vspace*{4mm}

This article is organized as follows.

In the Section \ref{sec:definition} we construct the model of random interlacements on transient weighted graphs and show that (\ref{eq:Q_u}) characterizes the image measure $Q^u$ of the interlacement set, see Remark~\ref{prop:Q_u}.

We prove in Section~\ref{sec:harry}, Theorem~\ref{th:harry} the Harris-FKG inequality for $Q^u$ and we state in Remarks~\ref{rem:harrisfinite} 1) and 2) the main obstructions concerning the use of standard techniques to prove this theorem.

In Section \ref{sec:nonamenable} we establish two results based on isoperimetric inequalities. Theorem~\ref{th:nonamenable} proves the claim (\ref{eq:statnonamenable}) whereas Theorem~\ref{th:isoper} shows (\ref{eq:statis6}).

In Section~\ref{sec:trees} we prove Theorem~\ref{th:tree} yielding (\ref{eq:stattree}) and also Corollary~\ref{cor:regular}, which exhibits the explicit formula (\ref{eq:criticaltreeintro}) for the critical value of regular trees.

Finally let us comment on our use of constants. Throughout this paper, $c$ will be used to denote a positive constant depending only on the graph under consideration, which can change from place to place. We write $c_1, \dots, c_5$ for fixed positive constants (also depending only on the graph under consideration), which refer to their first appearance in the text.

\textbf{Acknowledgments} - We are grateful to Alain-Sol Sznitman for the important corrections and encouragement.

\section{Definition of the model}
\label{sec:definition}

In this section we define the model of random interlacements in the more general setting of transient weighted graphs and prove the characterization of the process in terms of (\ref{eq:Q_u}).

We first introduce some further notation. We let $\lfloor \cdot \rfloor$ denote the integer part of a positive real number. Given two configurations $\alpha, \alpha' \in \{0,1\}^V$ we write $\alpha \succcurlyeq \alpha'$ if $\alpha(x) \geqslant \alpha'(x)$ for all $x \in V$. We say that a function $f$ defined on $\{0,1\}^V$ is increasing, if $f(\alpha) \geqslant f(\alpha')$, whenever $\alpha \succcurlyeq \alpha'$. We also denote, for real $a < b$, the uniform probability on the interval $[a,b]$ by $\mathcal{U}[a,b]$.

For a graph $\G = (V, \E)$, we say that $x, y \in V$ are neighbors (and write $x \leftrightarrow y$) if $\{x,y\} \in \E$. The \textit{degree} of a vertex $x \in V$ is defined as the number of edges incident to $x$.

Thoughout this article, the term \textit{path} always denotes nearest-neighbor finite paths, i.e. $\tau: \{0,\cdots,n\} \rightarrow V$ such that $\tau(i) \leftrightarrow \tau(i + 1)$ for $0 \leqslant i < n$, where $n \geqslant 0$ is what we call the \textit{length} of the path. We denote by $d_{\G}(x,y)$ (or simply $d(x,y)$ in case there is no ambiguity) the distance between $x$ and $y$, which is the smaller length among paths from $x$ to $y$. We write $B(x,n) = \{z \in V; d(z,x) \leqslant n\}$ for the closed ball with center $x \in V$ and radius $n$. For a set $K \subset V$, $d(x,K)$ stands for the distance between $x$ and $K$, i.e. the infimum of the distances between $x$ and the elements of $K$. The boundary of $K$, $\partial K$, is the set of points of $K$ that have some neighbor in $K^c$. The cardinality of $K$ is denoted by $|K|$.

Suppose that $\G$ is endowed with some weight function $\C$, see the definition in the introduction. For a finite set $A \subset V$, we define its capacity as
\begin{equation}
\label{eq:capacity}
\text{cap}(A) = \inf\bigg\{\frac{1}{2}\sum_{x,y \in V} |f(x) - f(y)|^2\C_{x,y}; f \text{ has finite support, } f \equiv 1 \text{ in } A\bigg\}.
\end{equation}
We call a graph transient when the capacity of some singleton (equivalently any finite set) is positive, see \cite{woess}, Theorem~2.12. From now on, we always assume that the weighted graph under consideration is transient.

The space $W_+$ stands for the set of infinite trajectories, that spend only a finite time in finite sets
\begin{equation}
\begin{split}
W_+ = \big\{ \gamma : \mathbb{N} \rightarrow V ; &\gamma (n) \leftrightarrow \gamma (n+1) \mbox{ for each } n \geqslant 0 \mbox{ and }\\
&\{n; \gamma(n) = y\} \text { is finite for all } y \in V \big\}.\\
\end{split}
\end{equation}
We endow $W_+$ with the $\sigma$-algebra $\mathcal{W}_+$ generated by the canonical coordinate maps $X_n$.
For $w \in W_+$ and $K \subset V$, we write $\widetilde{H}_K$ for the hitting time of $K$ by the trajectory $w$
\begin{equation}
\label{eq:hittime}
\widetilde{H}_K(w) = \inf \{n \geqslant 1; X_n(w) \in K\}.
\end{equation}

Given a weighted graph, we let $P_x$ stand for the law of the random walk associated with the transition matrix $q(\cdot,\cdot)$ starting at some point $x \in V$, see the definition in the introduction. If $\rho$ is some measure on $V$, we write $P_\rho = \sum_{x \in V} P_x \rho(x)$.

The assumed transience of the weighted graph $\G$, see below (\ref{eq:capacity}), is equivalent to the transience of the associated random walk. This means that $P_x[\widetilde{H}_{\{x\}} = \infty] > 0$ for every $x \in V$, we quote \cite{woess} Theorem~2.12. As a consequence, the set $W_+$ supports $P_x$.

If we define the \textit{equilibrium measure} by
\begin{equation}
\label{eq:equilibmeasure}
e_{A}(x) = 1_{\{x \in A\}}P_x[\widetilde{H}_A = \infty] \cdot \mu_x,
\end{equation}
we have the equality:
\begin{equation}
\label{eq:capeq}
\text{cap}(K) = \sum_{x \in K} e_K(x)
\end{equation}
for any finite set $K \subset V$, see for instance \cite{david2}, Proposition~2.3.

We further consider the space of doubly infinite trajectories that spend only a finite time in finite subsets of $V$
\begin{equation}
\begin{split}
W = \big\{ \gamma:\mathbb{Z} \rightarrow V ; &\gamma(n) \leftrightarrow \gamma(n+1) \mbox{ for each } n \in \mathbb{Z}, \mbox{ and }\\
&\{n; \gamma(i) = y\} \text { is finite for all } y \in V \big\}.
\end{split}
\end{equation}

We define the shift map, $\theta_k: W \rightarrow W$
\begin{equation}
\theta_k(w)(\cdot) = w(\cdot + k), \text{ for } k \in \mathbb{Z}.
\end{equation}
And for $w \in W$, we define the entrance time $H_K$ in a set $K \subset V$ by
\begin{equation}
H_K(w) = \inf \{n \in \mathbb{Z}; X_n(w) \in K\}.
\end{equation}

Consider the space $W^*$ of trajectories in $W$ modulo time shift
\begin{equation}
\label{eq:W_*}
W^* = W / \sim, \text{ where } w \sim w' \iff w(\cdot) = w'(\cdot + k) \text{ for some } k \in \mathbb{Z}. 
\end{equation}
and denote with $\pi^*$ the canonical projection from $W$ to $W^*$. The map $\pi^*$ induces a $\sigma$-algebra in $W^*$ given by $\mathcal{W}^* = \{A \subset W^*; (\pi^*)^{-1}(A) \in \mathcal{W}\}$, which is the largest $\sigma$-algebra on $W^*$ for which $(W,\mathcal{W}) \overset{\pi^*}{\rightarrow} (W^*,\mathcal{W}^*)$ is measurable.

Given a finite set $K \subset V$, write $W_K$ for the space of trajectories in $W$ that enter the set $K$, and denote with $W_K^*$ the image of $W_K$ under $\pi^*$.

The set of point measures on which one canonically defines the random interlacements is given by
\begin{equation}
\label{eq:omega}
\begin{split}
\Omega = \bigg\{\omega = \sum_{i\geqslant 1} \delta_{(w^*_i,u_i)}; & w^*_i \in W^*, u_i \in \mathbb{R}_+ \mbox{ and } \omega(W^*_K \times [0,u]) < \infty,\\
&\mbox{ for every finite } K \subset V \mbox{ and } u\geqslant 0 \bigg\}.
\end{split}
\end{equation}
It is endowed with the $\sigma$-algebra $\mathcal{A}$ generated by the evaluation maps $\omega \mapsto \omega(D)$
for $D \in$ \linebreak $\mathcal{W}^* \otimes \mathcal{B}(\mathcal{R}_+)$.

The interlacements are governed by a Poisson point process on $W^* \times \mathbb{R}_+$ with an intensity defined in terms of the measure $\nu$ on $W^*$, which is described in the theorem below. We also refer to \cite{weil} and \cite{silverstein} for measures similar to $\nu$, following the outline of \cite{hunt}. The construction we give here does not involve projective limits.

\begin{theorem}
\label{th:exist_nu}
There exists a unique $\sigma$-finite measure $\nu$ on $(W^*,\mathcal{W}^*)$ satisfying, for each finite set $K \subset V$,
\begin{equation}
\label{eq:nuQ}
1_{W^*_K} \cdot \nu = \pi^* \circ Q_K
\end{equation}
where the finite measure $Q_K$ on $W_K$ is determined by the following. Given $A$ and $B$ in $\mathcal{W}_+$ and a point $x \in V$,
\begin{equation}
\label{eq:Q_K}
Q_K[(X_{-n})_{n\geqslant 0} \in A, X_0 = x, (X_n)_{n\geqslant 0} \in B] = P_x[A|{\widetilde{H}}_K = \infty] e_K(x) P_x[B].
\end{equation}
\end{theorem}

\vspace{4 mm}

We follow the proof of Theorem~1.1 in \cite{sznitman} that establishes the existence of such measure in the case $V = \mathbb{Z}^d$.

\textit{Proof}. The uniqueness of $\nu$ satisfying (\ref{eq:nuQ}) is clear since, given a sequence of sets $K_n \uparrow V$, $W^* = \cup_n W^*_{K_n}$. For the existence, what we need to prove is that, for fixed $K \subset K' \subset V$,
\begin{equation}
\label{eq:Qcoerent}
\pi^* \circ (1_{W_K} \cdot Q_{K'}) = \pi^* \circ Q_K.
\end{equation}

We introduce the space
\begin{equation}
\label{eq:WKK}
W_{K,K'} = \{w \in W_K; H_{K'}(w) = 0\}
\end{equation}
and the bijection $s_{K,K'}:W_{K,K'} \rightarrow W_{K,K}$ given by 
\begin{equation}
\label{eq:s}
[s_{K,K'}(w)](\cdot) = w(H_K(w) + \cdot).
\end{equation}

To prove (\ref{eq:Qcoerent}), it is enough to show that
\begin{equation}
\label{eq:Qscoerent}
s \circ (1_{W_{K,K'}} \cdot Q_{K'}) = Q_K,
\end{equation}
where we wrote $s$ in place of $s_{K,K'}$ for simplicity. Indeed, by (\ref{eq:Q_K}), $1_{W_{K,K'}} \cdot Q_{K'} = 1_{W_K} \cdot Q_{K'}$ and one just applies $\pi^*$ on both sides of the equation above to obtain (\ref{eq:Qcoerent}).

We now consider the set $\Sigma$ of finite paths $\sigma:\{0,\cdots,N_\sigma\} \rightarrow V$, such that $\sigma(0) \in K'$, $\sigma(n) \notin K$ for $n < N_\sigma$ and $\sigma(N_\sigma)\in K$. We split the left hand-side of (\ref{eq:Qscoerent}) by partitioning $W_{K,K'}$ into the sets
\begin{equation}
\label{eq:WKs}
W_{K,K'}^\sigma = \{w \in W_{K,K'}; w \text{ restricted to } \{0,\cdots,N_\sigma\} \text{ equals } \sigma\}, \text{ for } \sigma \in \Sigma.
\end{equation}

For $w \in W_{K,K'}^\sigma$, we have $H_K(w) = N_\sigma$, so that we can write
\begin{equation}
\label{eq:splitQ}
s \circ (1_{W_{K,K'}} \cdot Q_{K'}) = \sum_{\sigma \in \Sigma} \theta_{N_\sigma} \circ (1_{W_{K,K'}^\sigma}\cdot Q_{K'}).
\end{equation}

To prove (\ref{eq:Qscoerent}), consider an arbitrary collection of sets $A_i \subset V$, for $i \in \mathbb{Z}$, such that $A_i \neq V$ for at most finitely many $i \in \mathbb{Z}$.
\begin{equation}
\begin{split}
s \circ (1_{W_{K,K'}} \cdot Q_{K'})[X_i \in A_i, i \in \mathbb{Z}] &=
\sum_{\sigma \in \Sigma} Q_{K'}[X_{i+N_\sigma}(w) \in A_i, i \in \mathbb{Z},w \in W_{K,K'}^\sigma] \\
&= \sum_{\sigma \in \Sigma} Q_{K'}[X_{i}(w) \in A_{i-N_\sigma}, i \in \mathbb{Z}, w \in W_{K,K'}^\sigma].
\end{split}
\end{equation}
Using the formula for $Q_{K'}$ given in (\ref{eq:Q_K}) and the Markov property, the above expression equals
\begin{equation}
\begin{aligned}
\label{eq:splitfinal}
& \textstyle\sum\limits_{x \in \text{Supp}(e_{K'})} \sum\limits_{\sigma \in \Sigma} \,\,
P_x[X_j \in A_{-j-N_\sigma}, j \geqslant 0, \widetilde{H}_{K'} = \infty] \cdot \mu_x \\
& \phantom{\sum\limits_{x \in \text{Supp}(e_{K'})}} \cdot P_x[X_n = \sigma(n) \in A_{n-N_\sigma}, 0 \leqslant n \leqslant N_\sigma] P_{\sigma(N_\sigma)}[X_n \in A_n, n \geqslant 0] \\
&= \textstyle\sum\limits_{\mbox{\fontsize{8}{8} \selectfont $\substack{y \in K\\ x \in \text{Supp}(e_{K'})}$}} 
\textstyle\sum\limits_{\sigma:\sigma(N_\sigma) = y}
P_x[X_j \in A_{-j-N_\sigma}, j \geqslant 0, \widetilde{H}_{K'} = \infty] \cdot \mu_x \\
&\phantom{\sum\limits_{x \in \text{Supp}(e_{K'})}} \qquad \cdot P_x[X_n = \sigma(n) \in A_{n-N_\sigma}, 0 \leqslant n \leqslant N_\sigma] P_{y}[X_n \in A_n, n \geqslant 0].
\end{aligned}
\end{equation}

For fixed $x \in \text{Supp}(e_{K'})$ and $y \in K$, we have
\begin{equation}
\label{eq:splity}
\begin{split}
& \textstyle\sum\limits_{\sigma:\sigma(N_\sigma) = y} P_x[X_j \in A_{-j-N_\sigma}, j \geqslant 0, \widetilde{H}_{K'} = \infty] \cdot \mu_x \\
& \phantom{\sum\limits_{x \in \text{Supp}(e_{K'})}} \cdot P_x[ X_n = \sigma(n) \in A_{n-N_\sigma},0\leqslant n \leqslant N_\sigma] \\
& \overset{\text{(reversibility)}}{=}\textstyle\sum\limits_{\mbox{\fontsize{8}{8} \selectfont $ \substack{\sigma:\sigma(N_\sigma) = y\\ \sigma(0) = x}$}}
P_x[X_j \in A_{-j-N_\sigma}, j \geqslant 0, \widetilde{H}_{K'} = \infty]\cdot \mu_y\\
&\phantom{\sum\limits_{x \in \text{Supp}(e_{K'})}} \qquad \quad \cdot P_y[X_m = \sigma(N_\sigma - m) \in A_{-m},0\leqslant m \leqslant N_\sigma]\\
&\overset{\text{(Markov)}}{=}\textstyle\sum\limits_{\mbox{\fontsize{8}{8} \selectfont $\substack{\sigma:\sigma(N_\sigma) = y\\ \sigma(0) = x}$}}
P_y[X_m = \sigma(N_\sigma - m) \in A_{-m},0\leqslant m \leqslant N_\sigma,\\
&\phantom{\sum\limits_{x \in \text{Supp}(e_{K'})}} \qquad \quad X_m \in A_{-m}, m \geqslant N_\sigma, \widetilde{H}_{K'} \circ \theta_{N_\sigma} = \infty]\cdot \mu_y\\
& \quad = P_y \Big[
\begin{array}{c}
\widetilde{H}_{K} = \infty, \text{ the last visit to $K'$}\\
\text{ occurs at $x$, } X_m \in A_{-m}, m \geqslant 0 \end{array} \Big] \cdot \mu_y.
\end{split}
\end{equation}
Using (\ref{eq:splity}) in (\ref{eq:splitfinal}) and summing over $x \in \text{Supp}(e_{K'})$, we obtain
\begin{equation}
\label{eq:insertsplit}
\begin{split}
s \circ (1_{W_{K,K'}} \cdot Q_{K'})[X_i \in A_i, i \in \mathbb{Z}] &= \textstyle\sum\limits_{y \in K}
P_y[\widetilde{H}_{K} = \infty, X_m = A_{-m}, m \geqslant 0] \cdot \mu_y\\
& \phantom{\textstyle\sum\limits_{x \in \text{Supp}(e_{K'})}}\cdot P_{y}[X_m \in A_m, m \geqslant 0]\\
&\overset{(\ref{eq:Q_K})}{=} Q_K[X_m \in A_m, m \in \mathbb{Z}].
\end{split}
\end{equation}
This shows (\ref{eq:Qscoerent}) and concludes the proof of the existence of the measure $\nu$ satisfying (\ref{eq:nuQ}). Moreover, $\nu$ is clearly $\sigma$-finite. $\square$

\vspace{4 mm}

\begin{remark}
\textnormal{
To recover the measure $\nu$ of \cite{sznitman} (as well as $\mathcal{I}^u$, see (\ref{eq:interlace}) below) one endows $\mathbb{Z}^d$ with the weight $(1/2d) \cdot 1_{x \leftrightarrow y}$, see (1.26) and Remark~1.4. Note that for this choice of conductance, $\mu$ is the counting measure on $\mathbb{Z}^d$.} $\square$
\end{remark}

\vspace{2 mm}

We are now ready to define the random interlacements. Consider on $\Omega \times \mathbb{R}_+$ the law $\mathbb{P}$ of a Poisson point process with intensity measure given by $\nu(dw^*) \otimes du$ (for a reference on this construction, see \cite{resnick} Proposition 3.6). We define the \textit{interlacement} and the \textit{vacant set} at level $u$ respectively as
\begin{gather}
\label{eq:interlace}
\mathcal{I}^u (\omega) = \bigg\{ \bigcup_{i; u_i \leqslant u} \mbox{Range}(w^*_i) \bigg\} \mbox{ and}\\
\label{eq:vacant}
\mathcal{V}^u(\omega) = V \setminus \mathcal{I}^u (\omega),
\end{gather}
for $\omega = \sum_{i\geqslant0} \delta_{(w^*_i,u_i)}$ in $\Omega$.

The next remark establishes the link with (\ref{eq:Q_u}).

\begin{remark}
\label{prop:Q_u}

\textnormal{1)
Consider, for $u \geqslant 0$, the map $\Pi^u:\Omega \rightarrow \{0,1\}^V$ given by
\begin{equation}
\label{eq:J}
(\Pi^u(\omega))_x = 1_{\{x \in \mathcal{V}^u(\omega)\}}, \text{ for $x$ in } V.
\end{equation}
Then the measure
\begin{equation}
\label{eq:QuPiP}
Q^u = \Pi^u \circ \mathbb{P}
\end{equation}
is characterized by (\ref{eq:Q_u}).}

\textnormal{
Indeed, for every finite subset $K$ of $V$ and $u \geqslant 0$, one has
\begin{equation}
\label{eq:remarkQu}
\begin{split}
Q^u[Y_x = 1 &\text{ for all } x \in K] = \mathbb{P}[\omega(W^*_K\times [0,u]) = 0] = \text{exp}(-u \cdot \nu(W^*_K))\\
&\overset{(\ref{eq:nuQ}), (\ref{eq:Q_K})}{=} \text{exp}\left(-u \cdot \textstyle{\sum_{x \in K}} e_K(x)\right) \overset{(\ref{eq:capeq})}{=} \text{exp}(-u \cdot \text{cap}(K)).
\end{split}
\end{equation}
Since the family of sets $[Y = 1 \text{ for all } x \in K]$ (where $K$ runs over all finite subsets of $V$) is closed under finite intersection and generates the $\sigma$-algebra $\mathcal{Y}$, (\ref{eq:Q_u}) uniquely determines $Q^u$.}

\textnormal{2)
Note also that in general the measure $Q^u$ neither dominates nor is dominated by any non-degenerate Bernoulli i.i.d. site percolation on $V$, see \cite{vladas}, Remark~1.1. However, as we will see in the next section, $Q^u$ satisfies the Harris-FKG inequality.} $\square$
\end{remark}

\section{The Harris-FKG inequality}
\label{sec:harry}

In this section we prove the Harris-FKG inequality (\ref{eq:harrisfkg}) for the measure $Q^u$, answering a question of \cite{sznitman} cf. Remark~1.6, 2).

A common strategy to prove (\ref{eq:harrisfkg}) is the following, see for instance \cite{grimmett} Theorem~2.4. In a first step, one proves that (\ref{eq:harrisfkg}) holds for random variables depending only on finitely many sites. A powerful sufficient condition in order to establish the inequality in the case of variables depending on finitely many coordinates is provided by the FKG-Theorem (see \cite{liggett} Corollary 2.12 p. 78). The general case is then handled by conditioning on the configuration inside some finite set $K$. The conditional expectations are then proved to be increasing random variables and one can apply the previous step. Finally one uses a martingale convergence argument to to achieve the result for the original random variables.

There are two main obstructions to using this strategy in our case. In Remark~\ref{rem:harrisfinite} 2), we prove that the sufficient condition (2.13) of \cite{liggett} p. 78 does not hold in general for the measures $Q^u$, so that the FKG-Theorem is not directly applicable. Further we provide in Remark~\ref{rem:harrisfinite} 1) an example in which the conditional expectation of an increasing function, with respect to the configuration of finitely many sites, is not increasing.

\vspace {4mm}

The strategy of the proof we present here strongly uses the construction of $Q^u$ in terms of the random interlacements. We consider an increasing sequence of finite sets $K_n \uparrow V$, which induces a certain sequence of $\sigma$-algebras $\mathcal{F}_{K_n}$. Roughly speaking, $\mathcal{F}_{K_n}$ keeps track of the behavior between the first and last visit to $K_n$ of all paths with level at most $u$ which meet the set $K_n$ . Given two function with finite second moment $f$ and $g$, we give an explicit representation of the conditional expectation of $f \circ \Pi^u$ and $g \circ \Pi^u$ with respect to $\mathcal{F}_{K_n}$ and prove that they are positively correlated. Finally, we prove (\ref{eq:harry}) by a martingale convergence argument.

The main theorem of this section comes in the following.

\begin{theorem}
\label{th:harry}
Consider $u \geqslant 0$ and let $f$ and $g$ be increasing random variables on $\{0,1\}^V$ with finite second moment with respect to the measure $Q^u$. Then one has
\begin{equation}
\label{eq:harry}
\int fg \, dQ^u \geqslant \int f \, dQ^u \int g \, dQ^u.
\end{equation}
\end{theorem}

\vspace{4 mm}

\textit{Proof.} We consider the countable space $\Gamma$ of finite paths in $V$. Given $K$ a finite subset of $V$, we define the functions $\phi_K: W^*_K \rightarrow \Gamma$ given by
\begin{equation}
\begin{split}
\phi_K (w^*) \text{ is the finite path starting when } & w^* \text{ first visits } K \text{ and}\\
\text{following $w^*$ step by step until its } &\text{last visit of $K$}.
\end{split}
\end{equation}
We also consider the partition of $W^*_K$ consisting of the sets
$$
W^*_{K,\gamma} = \{w^* \in W^*_{K}; \phi_K(w^*) = \gamma \}, \text{ for $\gamma \in \Gamma$}.
$$

Define the following random variables on $\Omega$
\begin{equation}
Z_{K,\gamma}(\omega)= \omega(W^*_{K,\gamma} \times [0,u]), \text{ for } \gamma \in \Gamma.
\end{equation}
They have Poisson distribution and are independent, since the sets $W^*_{K,\gamma}$, for $\gamma \in \Gamma$, are disjoint.

We regard $Z_K = (Z_{K,\gamma})_{\gamma \in \Gamma}$ as a random element of the space
$$
L = \{ (\alpha_\gamma)_{\gamma \in \Gamma} \in \mathbb{N}^\Gamma; \text{ $\alpha_\gamma$ is non-zero for finitely many $\gamma \in \Gamma$} \} \subset \mathbb{N}^\Gamma,
$$
with law denoted by $R_K$. Since the set $\mathbb{N}^{\Gamma}$ has a natural associated partial order, the notion of increasing and decreasing random variables on $L \subset \mathbb{N}^\Gamma$ is well defined.

Let $F$ and $G$ be two increasing random variables taking values in $L \subset \mathbb{N}^\Gamma$ that are square integrable with respect to $R_K$. We claim that
\begin{equation}
\label{eq:previousharri}
\int_L F G \, dR_K \geqslant \int_L F \, dR_K \int_L G \, dR_K.
\end{equation}
If $F$ and $G$ depend only on the value of finitely many coordinates $\Gamma' \subset \Gamma$, they can be trivially extended to increasing functions $F'$ and $G'$ in $\mathbb{N}^\Gamma$, one can choose for instance $F'(\alpha) = F(1_{\{\gamma \in \Gamma'\}} \cdot \alpha)$ and similarly for $G'$. In this case, since the law $R_K$ is a product measure on $\mathbb{N}^{\Gamma}$ (concentrated on $L$), (\ref{eq:previousharri}) is a direct consequence of the FKG-Theorem, see \cite{fortuin}, Proposition 1). The general case can be obtained by a martingale convergence argument as in \cite{grimmett} Theorem~2.4.

We define the $\sigma$-algebra $\mathcal{F}_K = \sigma(Z_K)$ and prove that it is possible to find increasing functions $F_K$ and $G_K$, on $L$, such that $F_K \circ Z_K = E[f \circ \Pi^u|\mathcal{F}_K]$ and $G_K \circ Z_K = E[g \circ \Pi^u|\mathcal{F}_K]$.

For this we construct the Poisson point process $\mathbb{P}$ in a more explicit way. Consider the probability measures on $W^*_{K,\gamma}$, defined by
\begin{equation}
\nu_{K,\gamma} = \frac{1_{W^*_{K,\gamma}} \cdot \nu}{\nu(W^*_{K,\gamma})}, \text{for $\gamma \in \Gamma$ such that $\nu(W^*_{K,\gamma}) > 0$, and arbitrarily otherwise}.
\end{equation}

On some auxiliar probability space $(S,\mathcal{S},\Sigma)$, we construct a collection of random elements $(\eta,\xi)_{\gamma,n}$ (for $\gamma \in \Gamma$ and $n \geqslant 0$) taking values on $W^*_{K,\gamma} \times [0,u]$. The law of this colection is characterized by the following:
\begin{equation}
\begin{array}{c}
\text{the $(\eta,\xi)_{\gamma\in \Gamma, n \geqslant 0}$ are independent and} \\
\text{each $(\eta,\xi)_{\gamma,n}$ is distributed as $\nu_{K,\gamma} \otimes \mathcal{U}[0,u]$}.\\
\end{array}
\end{equation}
In the same space $(S,\mathcal{S},\Sigma)$, independently of the collection above, we construct a Poisson point process, denoted by $N$, taking points on $(W^* \times \mathbb{R}_+) \setminus (W^*_K \times [0,u])$. More precisely, the process $N$ takes values in $\Omega' = \{ \omega \in \Omega; \omega(W^*_K \times [0,u]) = 0 \}$ and has intensity measure $1_{(W^*\times \mathbb{R}_+) \setminus (W^*_K \times [0,u])} \cdot \nu(dw^*)du$. Let $E^\Sigma$ denote the corresponding expectation.

Let $J_L$ and $J_S$ stand for the canonical projections on $L \times S$ and define the map:
\begin{equation}
\label{eq:defS}
\begin{split}
\Psi_K: L \times S & \longrightarrow \Omega \\
\left( \alpha, \sigma \right) & \longmapsto \sum_{\gamma \in \Gamma} \sum_{\,\, 0 < j \leqslant \alpha_\gamma} \delta_{(\eta, \xi)_{\gamma, j}(\sigma)} + N(\sigma).
\end{split}
\end{equation}
With these definitions, one has
\begin{equation}
\label{eq:Zandalpha}
J_L = Z_K \circ \Psi_K.
\end{equation}

It follows from the procedure to construct a Poisson point process, see for instance \cite{resnick} Proposition~3.6 p. 130, that
\begin{equation}
\label{eq:simulation}
\Psi_K \circ (R_K \otimes \Sigma) = \mathbb{P}.
\end{equation}

For a given $\alpha \in L$, choose $F_K(\alpha)$ and $G_K(\alpha)$ as $E^\Sigma [f \circ \Pi^u \circ \Psi_K(\alpha, \sigma)]$ and $E^\Sigma [g \circ \Pi^u \circ \Psi_K(\alpha, \sigma)]$ respectively. We now check that $F_K \circ Z_K$ and $G_K \circ Z_K$ are versions of the conditional expectations of $f \circ \Pi^u$ and $g \circ \Pi^u$ with respect to $\mathcal{F}_K$. Indeed, denoting with $E$ the expectation relative to $R_K \otimes \Sigma$, we find that for $\alpha_0 \in L$,
\begin{equation}
\begin{split}
\mathbb{E}[Z_K = \alpha_0, F_K \circ Z_K] & = \mathbb{P}[Z_K = \alpha_0] E^\Sigma[f \circ \Pi^u \circ \Psi_K(\alpha_0,\sigma)]\\
& \overset{(\ref{eq:Zandalpha}),(\ref{eq:simulation})}{=} E [J_L = \alpha_0] E[f \circ \Pi^u \circ \Psi_K(\alpha_0,J_S)]\\
& \overset{\text{(\ref{eq:Zandalpha}), indep.}}{=} E [Z_K \circ \Psi_K = \alpha_0, f \circ \Pi^u \circ \Psi_K]\\
& \overset{(\ref{eq:simulation})}{=} \mathbb{E} [Z_K = \alpha_0, f \circ \Pi^u].
\end{split}
\end{equation}

Note that if $\alpha, \alpha' \in L$ are such that $\alpha_\gamma \geqslant \alpha'_\gamma$ for every $\gamma \in \Gamma$, we have by (\ref{eq:defS}) that
\begin{equation}
\begin{split}
&\Pi^u \circ \Psi_K (\alpha, \cdot) \succcurlyeq \Pi^u \circ \Psi_K (\alpha', \cdot) \text{ for }  \text{every}\\
&\text{possible value of the second }  \text{coordinate}.
\end{split}
\end{equation}
Hence, the fact that $F_K$ and $G_K$ are increasing follows from the monotonicity of $f$ and $g$.

We now use (\ref{eq:previousharri}) and the fact that the conditional expectations of square integrable functions are also square integrable to deduce that
\begin{equation}
\begin{split}
\mathbb{E} \big[\mathbb{E}[f \circ \Pi^u |\mathcal{F}_K]\mathbb{E}[g \circ \Pi^u|\mathcal{F}_K] \big] & = \mathbb{E}\big[(F_K G_K) \circ Z_K\big]\\
& = \int F_K G_K \, dR_K\\
& \overset{(\ref{eq:previousharri})}{\geqslant} \int F_K \, dR_K \int G_K \, dR_K\\
& = \mathbb{E} \big[\mathbb{E}[f \circ \Pi^u|\mathcal{F}_K]\big] \mathbb{E} \big[\mathbb{E}[g \circ \Pi^u|\mathcal{F}_K]\big].
\end{split}
\end{equation}

We consider now a sequence $K_n \uparrow V$ and claim that $\mathcal{F}_{K_n}$ is an increasing sequence of $\sigma$-algebras. To see this, note that for $w^* \in W^*_{K_n}$, $\phi_{K_{n+1}}(w^*)$ determines $\phi_{K_n}(w^*)$. So that $Z_{K_{n+1}}$ contains all the necessary information to reconstruct $Z_{K_n}$.

Recall that for $x \in V$, $Y_x$ stands for the canonical coordinate on $\{0,1\}^V$. If $x \in K_n$, $Y_x \circ \Pi^u$ is determined by $Z_K$. As a result, $f \circ \Pi^u$ and $g \circ \Pi^u$ are both measurable with respect to $\sigma(\mathcal{F}_{K_n}; n \in \mathbb{N})$. The theorem now follows from (\ref{eq:QuPiP}) and the martingale convergence theorem. $\square$

\vspace{4 mm}

\begin{corollary}
\label{cor:uxu}
Given a transient, weighted graph $\G = (V,\E)$, the critical point $u_*$ in (\ref{eq:critical}) is well defined regardless of the choice of the base point $x$.
\end{corollary}

\textit{Proof.} Given $x, x' \in V$, as $\G$ is connected, we can choose a path $\tau$ joining $x$ to $x'$. Then we have:
\begin{equation}
\begin{split}
\eta(x,u) & = Q^u\big[\text{ the connected cluster of the set $\{z; Y_z = 1\}$ containing $x$ is infinite}\big]\\
& \geqslant Q^u\big[\text{ the connected cluster of the set $\{z; Y_z = 1\}$ containing $x'$ is infinite}\\
& \quad \quad \quad \text{ and contains } Range(\tau) \big]\\
& \overset{(\ref{eq:harry})}{\geqslant} \eta(x',u) \cdot Q^u\big[ Range(\tau) \subset \{z; Y_z = 1\}]\\
& \overset{(\ref{eq:Q_u})}{=} \eta(x',u) \cdot \text{exp}\{-u \cdot \text{Cap}(Range(\tau))\}.
\end{split}
\end{equation}
Hence, $\eta(x',u) > 0$ implies $\eta(x,u) > 0$. The same argument in the opposite direction thus shows that the positivity of $\eta(x,u)$ does not depend on the point $x$. $\square$

\vspace{4 mm}

\begin{remark}
\label{rem:harrisfinite}
\textnormal{1) As mentioned at the beginning of this section, given an increasing function on $\{0,1\}^V$, we cannot always choose an increasing version of its $Q^u$-conditional expectation with respect to the configuration of finitely many sites. One example is given in the following.}

\textnormal{Take $V$ as $\mathbb{N}$ and connect with an edge all the points of $\mathbb{N}$ within distance $1$ to each other and also $0$ to $3$ as in the Figure~\ref{fig:weighted2}. For $n \geqslant 3$, we assign the weight $e^n$ to the edge that joins $n$ and $n+1$, and weight $1$ for the rest of the edges. It is known that the random walk on $V$ induced by the conductions defined above is transient, hence we can define the interlacement process on it.}

\begin{figure} [ht]
\psfrag{a-1}{$0$}
\psfrag{a0}{$2$}
\psfrag{a1}{$3$}
\psfrag{a2}{$4$}
\psfrag{a00}{$1$}
\begin{center}
\includegraphics[angle=0, width=0.3\textwidth]{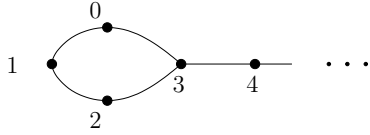}\\
\caption{The graph defined in the Remark~\ref{rem:harrisfinite}.}\label{fig:weighted2}
\end{center}
\end{figure}

\textnormal{We show that the events $[(Y_1,Y_2,Y_3) = (0,0,0)]$ and $[(Y_1,Y_2,Y_3) = (0,1,0)]$ have positive $Q^u$-probability. Moreover, despite of the fact that $Y_0$ is an increasing function on $\{0,1\}^\mathbb{N}$, we prove that}
\begin{equation}
\label{eq:expect}
E^u\big[Y_0 | (Y_1,Y_2,Y_3) = (0,0,0)\big] > 0 = E^u\big[Y_0 | (Y_1,Y_2,Y_3) = (0,1,0)\big].
\end{equation}

\textnormal{First consider the set}
\begin{equation}
W^*_r = \big\{w^* \in W^*; Range(w^*) \cap \{0,1,2,3\} = \{1,2,3\}\big\}.
\end{equation}
\textnormal{Which is disjoint from $W^*_{\{0\}}$. Now, using (\ref{eq:remarkQu}), we have}
\begin{equation}
\label{eq:0123}
\begin{split}
Q^u\big[(Y_0,Y_1,Y_2,Y_3) = (1,0,0,0)\big] \geqslant \, & \mathbb{P}\big[\omega(W^*_r \times [0,u]) > 0\big] \cdot \mathbb{P}\big[\omega(W^*_{\{0\}} \times [0,u]) = 0\big]\\
= \, &(1-e^{-u \cdot \nu(W^*_r)})e^{-u \cdot \nu(W^*_{\{0\}})},
\end{split}
\end{equation}
\textnormal{which is positive. Indeed, one easily checks that $W^*_r = W^*_{\{1\}} \setminus W^*_{\{0\}}$, and by (\ref{eq:nuQ}), (\ref{eq:Q_K}) we conclude that}
\begin{equation}
\label{eq:measureWr}
\begin{split}
\nu \left(W^*_{\{1\}} \setminus W^*_{\{0\}} \right) \geqslant P_1 [H_{\{0\}} = \infty] \cdot e_K(1) \cdot P_1[H_{\{0\}} = \infty | H_{\{1\}} = \infty ] > 0.
\end{split}
\end{equation}
\textnormal{This gives us the left-hand inequality of (\ref{eq:expect}).}

\textnormal{We claim that the configuration $(Y_1,Y_2,Y_3)=(0,1,0)$ has positive probability. Indeed, the exchange of the vertices $0$ and $2$ defines an isomorphism of the weighted graph, so that}
\begin{equation}
\label{eq:123}
\begin{split}
Q^u[(Y_1,Y_2,Y_3) = (0,1,0)] \geqslant & Q^u[(Y_0,Y_1,Y_2,Y_3) = (0,0,1,0)]\\
= & Q^u[(Y_0,Y_1,Y_2,Y_3) = (1,0,0,0)] \overset{(\ref{eq:0123}), (\ref{eq:measureWr})}{>} 0.
\end{split}
\end{equation}
\textnormal{On the other hand, the configuration $[(Y_0,Y_1,Y_2,Y_3) = (1,0,1,0)]$ is disjoint from the range of the map $\Pi^u$, since $\{1\}$ is a bounded component of $\{x; Y_x = 0\}$. Hence with (\ref{eq:remarkQu}) it has $Q^u$-probability zero, and the equality in (\ref{eq:expect}) follows.}

\vspace{4mm}

\textnormal{2) As we already mentioned, the Harris-FKG inequality for the measure $Q^u$ in general cannot be proved by the application of the FKG-Theorem (see \cite{liggett} Corollary 2.12 p. 78). Indeed, the condition (2.13) of \cite{liggett} p. 78}
\begin{equation}
\label{eq:conditpag78}
Q^u(\eta \wedge \zeta) Q^u(\eta \vee \zeta) > Q^u(\eta) Q^u(\zeta)
\end{equation}
\textnormal {does not hold in general for the measure $Q^u$. For instance, consider the example in Remark~\ref{rem:harrisfinite} 1) above. Define the two configurations for the state of the sites $(Y_0,Y_1,Y_2,Y_3)$}
$$
\eta = (0,0,1,0) \textnormal{ and } \zeta = (1,0,0,0).
$$
\textnormal{We know from (\ref{eq:0123}) and (\ref{eq:123}), that they have positive $Q^u$-probability. However the configuration $\eta \vee \zeta$ is given by $(1,0,1,0)$, and has zero $Q^u$-probability, contradicting (\ref{eq:conditpag78}).}

\vspace{4mm}

\textnormal{3) Finally we mention that the measure $Q^u$ does not satisfy the so-called Markov Field Property. This property states that, for every finite set $K \subset V$, the configuration on $(K \setminus \partial K)$ is independent of the configuration on the complement of $K$ when conditioned on what happens in $\partial K$.}

\textnormal{One can see from the example above, considering the set $K = \{1,2,3\}$, that}
\begin{equation}
Q^u[Y_2 = 1 | (Y_0, Y_1, Y_3) = (0,0,0)] \overset{(\ref{eq:123})}{>} 0 = Q^u[Y_2 = 1 | (Y_0, Y_1, Y_3) = (1,0,0)].
\end{equation}
\textnormal{This shows that the value of $Y_0$ can influence the value of $Y_2$, even if we condition on the configuration on $\partial K = \{1,3\}$. $\square$}

\end{remark}

\section{Some non-degeneracy results for $u_*$}
\label{sec:nonamenable}

In this section we derive some results on the non-degeneracy of the critical value $u_*$ under assumptions which involve isoperimetric inequalities. Roughly speaking, the main results of this section are Theorem~\ref{th:nonamenable} which shows the finiteness of $u_*$ for non-amenable graphs and Theorem~\ref{th:isoper} which shows the positivity of $u_*$ for $\G \times \mathbb{Z}$ when $\G$ satisfies $I \negmedspace S_6$.

We first introduce some further notation. Given two graphs $\G$ and $\G'$ with respective weight functions $\C$ and $\C'$, we define the weighted graph $\G \times \G'$ as follows. Two pairs $(x,x')$ and $(y,y')$ are considered to be neighbors if $x \leftrightarrow y$ and $x' = y'$ or if $x = y$ and $y \leftrightarrow y'$. In this case we define the weights between the two pairs by
\begin{equation}
\label{eq:productweight}
\C\left((x,x),(y,y')\right) = \C(x,y) 1_{\{x' = y'\}} + 1_{\{x = y\}}\C'(x',y').
\end{equation}

The first theorem of this section shows the finiteness of the critical value $u_*$ for graphs with bounded degrees, with weights bounded from above and from below and satisfying the strong isoperimetric inequality.

Loosely speaking, in the proof we first derive an exponential bound for the probability that a given set is vacant, see (\ref{eq:expbound}). With this bound it is straightforward to control the growth of the number of self avoiding paths of length $n$, starting from some fixed point. This proof does not apply for general transient weighted graphs, for instance in $\mathbb{Z}^d$ with $d \geqslant 3$, because a bound of type (\ref{eq:expbound}) fails in general. Indeed, according to (\ref{eq:Q_u}), $\mathbb{P}[A \subset \mathcal{V}^u]$ is given by $\text{exp} (-u \cdot \text{cap}(A))$ and the capacity of an arbitrary subset of $\mathbb{Z}^d$ is not bounded from below by any linear function of $|A|$, see \cite{sznitman} Remark~1.6 1).

\begin{theorem}
\label{th:nonamenable}
Let $\G$ be a graph of degree bounded by $q$, endowed with a weight function $\C$ which is bounded from below by $m$, satisfying the strong isoperimetric inequality. Then
\begin{equation}
u_* < \infty.
\end{equation}
In other words, for sufficiently large values of $u$, the vacant set $\mathcal{V}^u$ does not percolate.
\end{theorem}

\vspace{4 mm}

\textit{Proof.} We first note that every non-amenable graph is automatically transient. This follows for instance from (\ref{eq:boundcapacity}) below.

A weighted graph satisfies the strong isoperimetric inequality (i.e. (\ref{eq:strongperim}) with $d = \infty$) if and only if exists a $\bar k > 0$ such that the graph satisfies the Dirichlet inequality
\begin{equation}
\label{eq:sobineq}
\Vert f \Vert _2^2 \leqslant \bar \kappa \cdot \frac {1}{2} \sum_{x,y \in V} |f(x) - f(y)|^2\C_{x,y}, \text{ where $\Vert f \Vert _2 = \bigg( \sum_{x \in V} f(x) \mu_x \bigg)^{1/2}$,}
\end{equation}
for every $f:V \rightarrow \mathbb{R}$ with compact support. See for instance \cite{woess} Theorem~10.3.

From the definition of capacity (\ref{eq:capacity}), for any finite $A \subset V$,
\begin{equation}
\label{eq:boundcapacity}
\text{cap}(A) \geqslant \inf \bigg\{ \frac{1}{\bar \kappa} \cdot \Vert f \Vert _2^2; f \equiv 1 \text{ in } A, \text{ $f$ with compact support} \bigg\} \geqslant \frac{1}{\bar \kappa} \mu(A).
\end{equation}
Inserting this in the equation (\ref{eq:Q_u}) and using the lower bound on $\C$, we obtain, setting $\beta = \bar \kappa^{-1}$, the desired exponential bound
\begin{equation}
\label{eq:expon}
\mathbb{P}[A \subset \mathcal{V}^u] = Q^u[Y_x = 1, \text{ for all } x \in A] \leqslant \text{exp}(-u \cdot \beta \mu(A)) \leqslant \text{exp}(-u \cdot \beta \cdot m \cdot n).
\end{equation}

Given a fixed point $x_o \in V$, we can bound $\eta(u,x_o)$ (recall the definition in (\ref{eq:eta})) by the probability that $x_o$ is connected in $\mathcal{V}^u$ to $\partial B(x_o,n)$. Using the fact that the graph $\G$ has bounded degree, we bound the latter probability by summing over the set $\Gamma_n$ of all self-avoiding paths starting at $x_o$ with length $n$. More precisely, we have
\begin{equation}
\label{eq:expbound}
\begin{split}
\eta(u,x_o) &\leqslant \mathbb{P}[x_o \text{ is connected in $\mathcal{V}^u$ to } \partial B(x_o,n)]\\
&\leqslant \mathbb{P} [\text{there exists a $\gamma \in \Gamma_n$; } Range(\gamma) \subset \mathcal{V}^u]\\
&\leqslant \sum_{\gamma \in \Gamma_n} \mathbb{P} [Range(\gamma) \subset \mathcal{V}^u]\\
&\leqslant q^n \cdot \text{exp}(-\beta \cdot u \cdot m \cdot n)
\end{split}
\end{equation}
and this last expression tends to zero with $n$ for $u > (\beta m)^{-1}$. This readily implies that the probability that the vacant set at level $u > (\beta m)^{-1}$ contains an infinite component is zero. $\square$

\vspace*{4 mm}

As mentioned in the introduction, one important motivation for the introduction of the random interlacements in \cite{sznitman} has been the study of the disconnection time of a discrete cylinder, see for instance \cite{sznitman_cil1} and \cite{sznitman_cil3}. A natural way to generalize this kind of disconnection problem is to consider cylinders with more general basis, see \cite{sznitman_how}. This motivates our next result. It establishes that if some graph $\G$ satisfying $I \negmedspace S_6$ (see \ref{eq:strongperim}), has bounded degree and weights bounded both from above and from below, then the critical value $u_*$ for the graph $\G \times \mathbb{Z}$ is positive.

For the proof, we rely on some classical results of random walks on graphs to obtain a bound on the Green function of $\G \times \mathbb{Z}$. The rest of the proof is an adaptation of the renormalization argument in \cite{sznitman}, Proposition 4.1. 


\begin{theorem}
\label{th:isoper}
Let $\G$ be a graph of bounded degree endowed with weights that are bounded from above and from below, satisfying $I \negmedspace S_6$. Then the critical value $u_*$ of the graph $\G \times \mathbb{Z}$ is positive.
\end{theorem}

\textit{Proof.} Again, the transience of $\G \times \mathbb{Z}$ follows from the assumptions on $\G$, see the equation (\ref{eq:hittingbound}) below and the claim above (\ref{eq:equilibmeasure}).

Since $\G$ satisfies $I \negmedspace S_6$ and $\mathbb{Z}$ satisfies $I \negmedspace S_1$, the product $\G \times \mathbb{Z}$ satisfies $I \negmedspace S_7$, see \cite{woess}, 4.10 p. 44. This implies the upper bound
\begin{equation}
\label{eq:visitbound}
\sup_{x,y \in V} P_x[X_n = y] \leqslant C n^{-7/2},
\end{equation}
for some $C > 0$, see for instance \cite{woess}, 14.3 and 14.5 (a) p. 148. The bound (\ref{eq:visitbound}) implies the heat kernel bound
\begin{equation}
\label{eq:heatbound}
P_x[X_n = y] \leqslant C_1 n^{-7/2} \text{exp}\left\{ \frac{-d_{\G}(x,y)^2}{C_2n} \right\}, (x, y) \in \G \times \mathbb{Z},
\end{equation}
for some $C_1, C_2 > 0$, see \cite{woess}, 14.12 p. 153.

From this, by similar estimates as in \cite{lawler}, 1.5.4 p. 31, one obtains for some $C_3 > 0$ and all $(x,y) \in \G \times \mathbb{Z}$,
\begin{equation}
\label{eq:hittingbound}
P_x[H_{\{y\}} < \infty] \leqslant \sum_{n \geqslant 0} P_x[X_n = y] \leqslant \frac{C_3}{d_{\G}(x,y)^5}.
\end{equation}

We will use this bound in the renormalization argument we mentioned above. This renormalization will take place on an isometric copy of the upper plane $\mathbb{Z}_+ \times \mathbb{Z}$ that we find in $\G \times \mathbb{Z}$. More precisely, take a path $\tau: \mathbb{N} \rightarrow V$ satisfying what we call the \textit{half-axis property}:
\begin{equation}
\label{eq:taudistance}
d_{\G}(\tau(n), \tau(m)) = |n-m|, \text{ for all } n,m \in \mathbb{N}.
\end{equation}
The existence of such path is provided by \cite{watkins} Theorem 3.1.

We can now find an isometry between $\mathbb{Z}_+ \times \mathbb{Z}$ and a subset of $\G \times \mathbb{Z}$, with respect to the graph distances $d_{\, \mathbb{Z}_+ \times \mathbb{Z}}(\cdot,\cdot)$ and $d_{\, \G \times \mathbb{Z}}(\cdot,\cdot)$ which are defined as above (\ref{eq:capacity}). We take the map that, for a given pair $(i,j) \in \mathbb{Z}_+ \times \mathbb{Z}$, associates $(\tau(i),j)$.

To see why this defines an isometry, we use (\ref{eq:taudistance}) and note that for two graphs $\G_1$ and $\G_2$, one has $d_{\, \G_1 \times \G_2}((i,j),(i',j')) = d_{\, \G_1}(i,i') + d_{\, \G_2}(j,j')$. Indeed, if one concatenates two minimal paths, the first joining $(i,j)$ to $(i',j)$ in $\G_1 \times \{j\}$ (which is an isometric copy of $\G_1$) and the second joining $(i',j)$ to $(i',j')$ in $\{i'\} \times \G_2$ (which is isometric to $\G_2$), one obtains $d_{\, \G_1 \times \G_2}((i,j),(i',j')) \leqslant d_{\, \G_1}(i,i') + d_{\, \G_2}(j,j')$. To prove the other inequality, we note that for every path joining $(i,j)$ and $(i',j')$ one can decompose it into its horizontal and vertical steps (corresponding respectively to steps between pairs $(k,l) \leftrightarrow (k',l)$ and $(k,l) \leftrightarrow (k,l')$, see above (\ref{eq:productweight})) to obtain paths in $\G_1$ and $\G_2$ joining $i$ to $i'$ and $j$ to $j'$, respectively.

From now on we make no distinction between $\mathbb{Z}_+ \times \mathbb{Z}$ and $S = Range(\tau) \times \mathbb{Z} \subset V \times \mathbb{Z}$.

We say that $\tau:\{0,\cdots,n\} \rightarrow \mathbb{Z}_+ \times \mathbb{Z}$ is a $*$-path if
$$
|\tau(k + 1) - \tau (k)|_\infty = 1, \text{ for all } k \in \{0,\cdots, n-1\},
$$
where $|p|_\infty$ is the maximum of the absolute value of the two coordinates of $p \in \mathbb{Z}^2$

The rest of the proof follows the argument for the Proposition 4.1 in \cite{sznitman} with some minor modifications. For the reader's convenience, we write here the proof together with the necessary adaptions. Define
\begin{equation}
\label{eq:L_n}
L_0 > 1, L_{n+1} = l_nL_n, \text{ where } l_n = 100 \lfloor L^a_n \rfloor \text{ and } a = \frac{1}{1000}.
\end{equation}
We will consider a sequence of boxes in $S$ of size $L_n$ and define the set of indexes
\begin{equation}
\label{eq:indexJn}
J_n = \{n\} \times (\mathbb{Z}_+ \times \mathbb{Z}), J_n' = \{n\} \times \left[(\mathbb{Z}_+ \setminus \{0\}) \times \mathbb{Z}\right].
\end{equation}
For $\overline{m} = (n,q) \in J_n$, we consider the box
\begin{equation}
\label{eq:box}
D_{\overline{m}} = (L_n q + [0, L_n)^2) \cap \mathbb{Z}^2,
\end{equation}
recalling that we have identified $S$ with $\mathbb{Z}_+ \times \mathbb{Z}$. And for $m \in J_n'$ we set
\begin{equation}
\label{eq:box2}
\widetilde{D}_m = \bigcup_{i,j \in \{-1,0,1\}} D_{(n,q + (i,j))}.
\end{equation}

Roughly speaking, our strategy is to prove that the probability of finding a $*$-path in the set $\mathcal{I}^u \cap S$ that separates the origin from infinite in $S$ is smaller than one. We do this by bounding the probabilities of the following crossing events
\begin{equation}
\label{eq:crossing}
B^u_{m} = \{\omega \in \Omega; \text{there exists in $I^u \cap S$ a $*$-path between $D_m$ and $\widetilde{D}_m$}\},
\end{equation}
where $m \in J_n'$. For $u > 0$ and $m \in J_n'$, we write
\begin{gather}
\label{eq:qum}
q^u_m = \mathbb{P}[B^u_n] \text{ and }\\
\label{eq:qun}
q^u_n = \sup_{\overline{m} \in J_n} q^u_{\overline{m}}.
\end{gather}

\begin{figure}
\psfrag{a}{$D_m$}
\psfrag{b}{$\widetilde{D}_m$}
\psfrag{c}{$\widetilde{D}_{\overline{m}_1}$}
\psfrag{d}{$D_{\overline{m}_1}$}
\psfrag{e}{$D_{\overline{m}_2}$}
\psfrag{f}{$\widetilde{D}_{\overline{m}_2}$}
\begin{center}
\includegraphics[angle=0, width=0.5\textwidth]{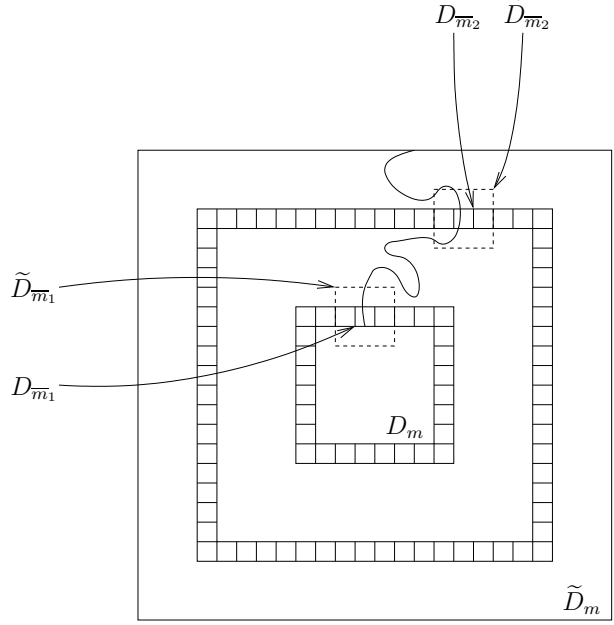}\\
\caption{The figure shows all the boxes with indexes in $\mathcal{K}_1$ and $\mathcal{K}_2$. Note that the event $B^u_m$ implies $B^u_{\overline{m}_1}$ and $B^u_{\overline{m}_2}$ for some $\overline{m}_1 \in \mathcal{K}_1$ and $\overline{m}_2 \in \mathcal{K}_2$.}\label{fig:weighted}
\end{center}
\end{figure}

In order to obtain an induction relation between $q^u_n$ and $q^u_{n+1}$ (that were defined in terms of two different scales) we consider for a fixed $m \in J_{n+1}'$ the indexes of boxes in the scale $n$ that are in the ``boundary of $D_m$''
\begin{equation}
\mathcal{K}^m_1 = \{\overline{m}_1 \in J_n; D_{\overline{m}_1} \subset D_m \text{ and $D_{\overline{m}_1}$ is neighbor of } S \setminus D_m\}.
\end{equation}
And the indexes of boxes at the scale $n$ and that have some point at distance $L_{n+1}/2$ of $D_m$
\begin{equation}
\mathcal{K}^m_2 = \{\overline{m}_2 \in J_n; D_{\overline{m}_2} \cap \mbox{\fontsize{10}{10} \selectfont $\{x \in S; d_{\G \times \mathbb{Z}}(z, D_m) = L_{n+1}/2\}$} \neq \emptyset \}.
\end{equation}
The boxes associated with the two sets of indexes above are shown in Figure~\ref{fig:weighted}. In this figure we also illustrate that the event $B^u_m$ implies the occurrence of both $B^u_{\overline{m}_1}$ and $B^u_{\overline{m}_2}$ for some choice of $\overline{m}_1 \in \mathcal{K}^m_1$ and $\overline{m}_2 \in \mathcal{K}^m_2$.

This, with a rough counting argument, allows us to conclude that
\begin{equation}
\label{eq:boundq}
q^u_m \leqslant cl^2_n \sup_{{\mbox{\fontsize{8}{8} \selectfont $\substack{\overline{m}_1 \in \mathcal{K}^m_1\\ \overline{m}_2 \in \mathcal{K}^m_2}$}}} \mathbb{P}[B^u_{\overline{m}_1} \cap B^u_{\overline{m}_2}], \text{ for all } u \geqslant 0.
\end{equation}
We now want to control the dependence of the process in the two boxes $\widetilde{D}_{\overline{m}_1}$ and $\widetilde{D}_{\overline{m}_2}$. For this we will use the estimates (\ref{eq:hittingbound}) and split the set $W^*$ as follows.

Let $V$ denote $\widetilde{D}_{\overline{m}_1} \cup \widetilde{D}_{\overline{m}_2}$ and write the set $W^*_V$ as the union below
\begin{equation}
\label{eq:splitW_V}
W^*_V = W^*_{1,1} \cup W^*_{1,2} \cup W^*_{2,1} \cup W^*_{2,2},
\end{equation}
where
\begin{equation}
\label{eq:W_ij}
\begin{split}
W^*_{1,1} &= W^*_{\widetilde{D}_{\overline{m}_1}} \setminus W^*_{\widetilde{D}_{\overline{m}_2}}\\
W^*_{1,2} = \{\omega \in W^*_{\widetilde{D}_{\overline{m}_1}}&; X_{H_V}(\omega) \in \widetilde{D}_{\overline{m}_1} \text{ and } H_{\widetilde{D}_{\overline{m}_2}} < \infty\}.
\end{split}
\end{equation}
And $W_{2,2}$, $W_{2,1}$ defined with similar formulas. Note that the union above is disjoint, so that
\begin{equation}
\label{eq:poissindep}
\text{the Poisson point processes obtained by restricting $\omega$ to $W^*_{i,j}$ are independent}.
\end{equation}

For some measurable $A \subset W^*$ and some index $m \in J_{n+1}'$ with $n \geqslant 0$ we introduce the notation
\begin{equation}
\label{eq:Bum_delta}
B^u_m(A) = \{\omega \in \Omega; 1_{A \times [0,u]} \cdot \omega \in B^u_m\}
\end{equation}

Since the trajectories in the set $W^*_{i,i}$ are disjoint of the box $\overline{m}_j$ for $i,j \in \{0,1\}$ and $i \neq j$, we conclude that
\begin{equation}
\label{eq:BeqBdelta}
B^u_{\overline{m}_1} = B^u_{\overline{m}_1}\left(\cup_{i,j \in \{0,1\}} W_{i,j} \right) = B^u_{\overline{m}_1}\left(W^*_{1,1} \cup W^*_{1,2} \cup W^*_{2,1}\right), \text{ for } u \geqslant 0
\end{equation}
and a similar formula for $B^u_{\overline{m}_2}$.

We now control the dependence of the process in the two different boxes. This is done by simply requiring that no trajectory reaches both $D_{\overline{m}_1}$ and $D_{\overline{m}_2}$.
\begin{equation}
\label{eq:bounddependence}
\begin{split}
\mathbb{P}[B^u_{\overline{m}_1} \cap B^u_{\overline{m}_1}] &= \mathbb{P}\left[B^u_{\overline{m}_1}\left(W^*_{1,1} \cup W^*_{1,2} \cup W^*_{2,1}\right) \cap B^u_{\overline{m}_2}\left(W^*_{2,2} \cup W^*_{2,1} \cup W^*_{1,2}\right)\right]\\
&\leqslant \mathbb{P}\left[B^u_{\overline{m}_1}\left(W^*_{1,1}\right) \cap B^u_{\overline{m}_2}\left(W^*_{2,2}\right), \omega(W_{1,2} \times [0,u]) = \omega(W_{2,1} \times [0,u]) = 0 \right]\\
& \quad + \mathbb{P}\left[\omega(W_{1,2} \times [0,u]) \text{ or } \omega(W_{2,1} \times [0,u]) \neq 0\right]\\
&\overset{(\ref{eq:poissindep})}{\leqslant} \mathbb{P}\left[B^u_{\overline{m}_1}\left(W^*_{1,1}\right)\right] \mathbb{P}\left[ B^u_{\overline{m}_2}\left(W^*_{2,2}\right) \right] + \mathbb{P}\left[\omega(W_{1,2} \times [0,u]) \neq 0\right]\\
& \quad + \mathbb{P}\left[\omega(W_{2,1} \times [0,u]) \neq 0\right]\\
&\leqslant q^u_{\overline{m}_1} q^u_{\overline{m}_2} + (1-\text{e}^{-u\cdot\nu(W_{1,2})}) + (1-\text{e}^{-u \cdot \nu(W_{2,1})})\\
&\overset{(\ref{eq:Q_K})}{\leqslant} q^u_{\overline{m}_1} q^u_{\overline{m}_2} + u\Big(P_{e_V}[X_0 \in \widetilde{D}_{\overline{m}_1}, H_{\widetilde{D}_{\overline{m}_2}} < \infty]\\
&\qquad \qquad \qquad \qquad + P_{e_V}[X_0 \in \widetilde{D}_{\overline{m}_2}, H_{\widetilde{D}_{\overline{m}_1}} < \infty]\Big)\\
&\overset{(\ref{eq:hittingbound})}{\leqslant} (q^u_n)^2 + cL^2_n \frac{L^2_n}{L^5_{n+1}}
\end{split}
\end{equation}
where we assumed in the last step that $u \leqslant 1$. Using (\ref{eq:boundq}) and taking the supremum over $m \in J_{n+1}$, we conclude that
\begin{equation}
\label{eq:boundq2}
q^u_{n+1} \leqslant c_1 l^2_n \left( (q^u_n)^2 + L^4_n L^{-5}_{n+1} \right).
\end{equation}

With help of this equation, we prove the next Lemma, which shows that for some choice of $L_0$ and $u$ taken small enough, $q^u_n$ goes to zero sufficiently fast with $n$.

\begin{lemma}
\label{lem:quntozero}
There exist $L_0$ and $\bar u = \bar u (L_0) > 0$, such that
\begin{equation}
\label{eq:qundecay}
q^u_n \leqslant \frac{c_3}{l^2_n L^{1/2}_n}
\end{equation}
for every $u < \bar u$.
\end{lemma}

\textit{Proof of the Lemma.} We define the sequence
\begin{equation}
\label{eq:defbn}
b_n = c_3 l^2_nq^u_n, \text{ for } n\geqslant 0.
\end{equation}
The equation (\ref{eq:boundq2}) can now be rewritten as
\begin{equation}
\label{eq:boundbn}
b_{n + 1} \leqslant c_4\left( \left( \frac{l_{n+1}}{l_n} \right)^2 b^2_n + (l_{n+1} l_n)^2 L^4_n L^{-5}_{n+1} \right), \text{ for } n \geqslant 0.
\end{equation}
With (\ref{eq:L_n}) one concludes that $(l_{n+1} l_n)^2 \leqslant cL^{2a}_nL^{2a}_{n+1} \leqslant cL^{4a + 2a^2}_n$. Inserting this in (\ref{eq:boundbn}) and using again (\ref{eq:L_n}), we obtain
\begin{equation}
b_{n+1} \leqslant c_5(L^{2a^2}_n b^2_n + L^{2a^2 - a - 1}_n) \leqslant c_5 L^{2a^2}_n(b^2_n + L^{-1}_n).
\end{equation}

We use this to show that, if for some $L_0 > (2c_5)^4$ and $u \leqslant 1$ we have $b_n \leqslant L^{-1/2}_n$, then the same inequality also holds for $n + 1$. Indeed, supposing $b_n \leqslant L^{-1/2}_n$, we have
\begin{equation}
\label{eq:inductionbn}
b_{n+1} \leqslant 2 c_5 L^{2a^2 - 1}_n \overset{(\ref{eq:L_n})}{\leqslant} 2c_5 L^{-1/2}_{n+1}L^{1/2(1+a) + 2a^2 - 1}_n \overset{(\ref{eq:L_n})}{\leqslant} 2c_5 L^{-1/2}_{n+1} L^{-1/4}_0 \leqslant L^{-1/2}_{n+1}.
\end{equation}
Which is the statement of the lemma. So all we still have to prove is that $b_0 \leqslant L^{-1/2}_0$ for $L_0 > (2c_5)^4$ and small enough $u$. Indeed,
\begin{equation}
\begin{split}
b_0 &\overset{(\ref{eq:defbn})}{=} c_3 l^2_0 q^u_0 \leqslant c_3 l^2_0 \sup_{m \in J_0} \mathbb{P}[\mathcal{I}^u \cap \widetilde{D}_m \neq \emptyset]\\
&\leqslant c_5 L^{2a + 2}_0 \sup_{x \in V} \mathbb{P}[x \in \mathcal{I}^u] \overset{(\ref{eq:Q_u}),(\ref{eq:capeq})}{\leqslant} c_5 L^{2a + 2}_0 (1 - e^{-\alpha u}),
\end{split}
\end{equation}
where $\alpha$ can be taken as $\sup_{x \in V} \mu_x$ since both the degree and the weights are bounded from above. For some $L_0 > (2c_5)^4$, we take $u(L_0)$ small enough such that $b_0 \leqslant L^{-1/2}_0$ for any $u \leqslant u(L_0)$. $\square$

\vspace*{4mm}

We now use this lemma to show that with positive probability, one can find an infinite connection from $(0,0)$ to infinite in the set $\mathcal{V}^u \cap S$. For this we choose $L_0$ and $u < u(L_0)$ as in the lemma. Writing $B_M$ for the set $[0,M] \times [-M,M] \subset S$, we have
\begin{equation}
\begin{split}
1 - \eta(u,(0,0)) &\leqslant \mathbb{P}[ (0,0) \text{ is not in an infinite component of } \mathcal{V}^u \cap S]\\
&\leqslant \mathbb{P}[\mathcal{I}^u \cap B_M \neq \emptyset] +
\begin{split}
\mathbb{P}\big[ &\text{there is a $*$-path in } S \setminus B_M\\
&\text{surrounding the point } (0,0) \text{ in } S \big]\\
\end{split}\\
&\leqslant \big(1 - \text{exp}(-u \cdot \text{cap}(B_M))\big)\\
&\quad + \sum_{n \geqslant n_0}
\begin{split}
\mathbb{P} \big[& \mathcal{I}^u \cap S \setminus B_M \text{ contains a $*$-path surrounding $(0,0)$ and}\\
& \text{passing through some point in $[L_n,L_{n+1} - 1] \times \{0\} \in S$} \big]
\end{split}
\end{split}
\end{equation}
The last sum can be bounded by $\sum_{n \geqslant n_0} \sum_m \mathbb{P}[B^u_m]$ where the index $m$ runs over all labels of boxes $D_m$ at level $n$ that intersect $[L_n,L_{n+1} - 1] \times \{0\} \subset S$. Since the number of such $m$'s is at most $l_n \leqslant c L^a_n$,
\begin{equation}
1 - \eta(u,(0,0)) \leqslant c L^2_{n_0}u + \sum_{n \geqslant n_0} cL^a_n L^{-1/2}_n \overset{(\ref{eq:L_n})}{\leqslant} c(L^2_{n_0}u + \sum_{n \geqslant n_0} L^{-1/4}_n).
\end{equation}
Choosing $n_0$ large and $u \leqslant u(L_0,n_0)$, we obtain that the percolation probability is positive. So that $u_* > 0$. $\square$



\section{Trees}
\label{sec:trees}

The next theorem roughly states that if $\G$ is a tree, the vacant cluster at level $u$ containing some given vertex $x \in V$ has the same law as the vacant cluster at $x$ under a Bernoulli independent site percolation with a suitable choice of occupancy probabilities (that depends on $x$ and $u$). A consequence of this fact is that $0 < u_* < \infty$, under some additional bounds on the degrees and weights of the tree.

The strategy of the proof is to partition $W^*$ into sets $W^{*,z}$ ($z \in V$) (see (\ref{eq:W_z}) below) and use $\omega(W^{*,z} \times [0,u])$ to build independent random variables that will induce the Bernoulli process. The main task in the proof is to check that the component containing $x$ is indeed the same both for the interlacement set and for the induced Bernoulli site percolation.

\begin{theorem}
\label{th:tree}
Let $\G = (V,\E)$ be a tree with bounded degree, endowed with a weight function which makes it a transient weighted graph and take $x \in V$ a fixed site. Consider the function $f_x :V \rightarrow [0,1]$ given by
\begin{equation}
\label{eq:bernoullicouple}
\begin{split}
f_x(z) = \thinspace &P_z\left[d(X_n,x) > d(z,x) \textnormal{ for every } n > 0 \right] \cdot \mu_z\\
&\cdot P_z\left[d(X_n,x) \geqslant d(z,x) \textnormal{ for every } n \geqslant 0 \right].
\end{split}
\end{equation}

Then the cluster of $\mathcal{V}^u$ containing $x$ in the interlacement process has the same law as the open cluster containing $x$ for the Bernoulli independent site percolation on $V$ characterized by $P[z \textnormal{ is open} \thinspace]$ $=$ $\textnormal{exp}(-u \cdot f_x(z))$.

In particular, if in addition the degree of each vertex of the tree is larger than or equal to three and the weights are bounded from above and from below, then the critical value $u_*$ is non-degenerate, i.e.
\begin{equation}
\label{eq:nondegofud}
0 < u_* < \infty.
\end{equation}
\end{theorem}

\vspace{4 mm}

\textit{Proof.} We introduce, the partition $W^* = \sqcup_{z \in V} W^{*,z}$, where
\begin{equation}
\label{eq:W_z}
W^{*,z} = \{w^* \in W^*_{\{z\}}; d(x,Range(w^*)) \geqslant d(x,z)\}.
\end{equation}

We claim that the sets $W^{*,z}$ are disjoint. Indeed, consider some trajectory $w^* \in W^{*,z}$. Since $\G$ is a tree, $z$ is the unique point in $Range(w^*)$ for which $d(x,z)$ attains the distance $d(x,Range(w^*))$. As a consequence, the random variables $(\omega(W^{*,z}\times [0,u]))_{z \in V}$ are independent. We use them to define a Bernoulli independent process given by
\begin{equation}
\label{eq:bernoulliproc}
Y^u_z(\omega) = 1_{[\omega(W^{*,z}\times[0,u]) \geqslant 1]} \text{ for } z \in V.
\end{equation}
Note by (\ref{eq:Q_K}) and (\ref{eq:W_z}) that the probability $\mathbb{P}[Y^u_z = 0]$ equals $\text{exp}(-u \cdot f_x(z))$.

We still need to prove that for the site percolation attached to the values of $Y^u_z$, the null-cluster containing $x$ coincides with the connected component of $\mathcal{V}^u$ containing $x$. From now on, we fix a point measure $\omega \in \Omega$ and write $Y^u_z$ and $\mathcal{I}^u$ instead of $Y^u_z(\omega)$ and $\mathcal{I}^u(\omega)$ for simplicity. With no loss of generality we assume $x \in \mathcal{V}^u$. Our claim will follow once we show that given a nearest neighbor self-avoiding path $\tau:\{0,\cdots,l\} \rightarrow V$ starting at $x$, $Range(\tau) \subset \mathcal{V}^u$ if and only if $Y_z = 0$ for every $z \in Range(\tau)$.

Suppose first that $Y_z = 0$ for every $z \in Range(\tau)$ and, by contradiction, that $\tau$ intersects $\mathcal{I}^u$. We call $\bar n \in \{0,\dots,l\}$ the first time at which $\tau$ meets $\mathcal{I}^u$. Since $1_{[x \in \mathcal{I}^u]} = 1_{[\omega(W^*_{\{x\}}\times [0,u]) \geqslant 1]} = Y_x = 0$, we conclude that $\bar n > 0$. Thus, $\tau(\bar n) \in \mathcal{I}^u$, but $\omega(W^{*,\tau(\bar n)} \times [0,u]) = 0$. Hence some $(w^*_i,u_i)$ with $u_i < u$ in the support of $\omega$ belongs to $W_{\{\tau(\bar n)\}}$ and gets closer to $x$ than $\tau(\bar n)$. Since $\G$ is a tree, the only fashion in which a trajectory of $W^*_{\{\tau(\bar n)\}}$ can get closer to $x$ than $\tau(\bar n)$ is by visiting the point $\tau(\bar n-1)$, which we know belongs to the vacant set, a contradiction.

Conversely, $Range(\tau) \subset \mathcal{V}^u$ implies that $1_{[\omega(W^*_{\{z\}} \times [0,u]) \geqslant 1]} = 0$ for every $z \in Range(\tau)$. Since $W^{*,z} \subset W^*_{\{z\}}$, we have $Y_z = 0$ for every $z \in Range(\tau)$. This concludes the proof of the first part of the theorem.

Let us now suppose in addition that every vertex of the tree has degree larger than or equal to three and that the weights are bounded from above and from below. We now prove the non-degeneracy of $u_*$ (\ref{eq:nondegofud}).

Since the tree has bounded degree, $f_x(z)$ is also bounded from above, say by $L$. And the construction above shows that the vacant cluster containing $x$ has the same law as under an independent Bernoulli site percolation on the tree. We can explore the cluster which contains $x$ in this process, in such a way that it corresponds to a Galton-Watson model for population growth. Every site $z$ in this cluster is regarded as an individual in the population, and the distance between $z$ and the progenitor $x$ is understood as the number of generations. The neighbors of $z$, in the next generation which are vacant in the Bernoulli process, are regarded as the descendents of $z$.

Since the degree of the tree is supposed to be larger than or equal to three, the distributions of descendants in this branching process stochastically dominates the sum of two independent Bernoulli variables with success probability $\text{exp}(-uL)$ (recall that $P[z \text{ is open} \thinspace]$ $=$ $\text{exp}\{(-u \cdot f_x(z)\} \geqslant \text{exp}\{-u \cdot L\}$). It follows from \cite{grimmett2}, 5.4 Theorem (5), that with positive probability the population of this branching process does not become extinct when $u < \text{log}(2) / L$. As a consequence, for $u < \text{log}(2) / L$, there is a positive probability that $x$ percolates in $\mathcal{V}^u$. This gives us the positivity of $u_*$. Here we did not use that the weights are bounded from below.

The finiteness of the critical value will follow from the Theorem~\ref{th:nonamenable} since $\G$ is non-amenable. Indeed, one can use the Theorem~(10.9) of \cite{woess} p. 114, to show that $\G$ endowed with the canonical weights is non-amenable. And the same will hold for $\G$ endowed with the original weights since they are bounded from above and from below and the degrees are bounded. $\square$


\vspace{4 mm}

\begin{corollary}
\label{cor:regular}
Let $\G$ be the regular tree with degree $d > 2$ endowed with weights $1/d$ for every edge. Then the critical point of the interlacement percolation is given by
\begin{equation}
\label{eq:criticaltree}
u_* = \frac{d(d-1)\text{log}(d-1)}{(d-2)^2}.
\end{equation}
\end{corollary}

\vspace{4 mm}

\textit{Proof}. We consider the construction provided by Theorem~\ref{th:tree} and give an explicit formula for $f_x(z)$. Given $z, z' \in V$ which are neighbors, one can use the strong Markov property on $H_{\{z\}}$ and the Lemma (1.24) from \cite{woess} p. 9, to conclude that
\begin{equation}
\label{eq:pz1z2}
P_{z}[H_{\{z'\}} < \infty] = \frac{\sum_{n \geqslant 0} P_{z}[X_n = z']}{\sum_{n \geqslant 0} P_{z}[X_n = z]} = \frac{d - \sqrt[]{d^2-4d+4}}{2(d-1)} = \frac{1}{d-1}.
\end{equation}

Now suppose that $z \neq x$ and the point $z'$ is the unique neighbor of $z$ for which the distance to $x$ equals $d(z, x) - 1$. From (\ref{eq:bernoullicouple}) and (\ref{eq:pz1z2}) one has
\begin{equation}
\begin{split}
f_x(z) & = P_z[X_1 \neq z', P_{X_1}[H_{\{z\}} = \infty]] \cdot \mu_z \cdot P_z[H_{\{z'\}} = \infty]\\
& =  \frac{d-1}{d} \cdot \frac{d-2}{d-1} \cdot \frac{d-2}{d-1} = \frac{(d-2)^2}{d(d-1)}.
\end{split}
\end{equation}

From this we conclude that the cluster containing $x$ has the same law as a branching process for which (except for the first generation) the numbers of descendants are binomial, i.i.d. random variables with expectation given by
\begin{equation}
\label{eq:expectationregular}
(d-1)\text{exp}\left(-u \frac{(d-2)^2}{d(d-1)}\right).
\end{equation}
The critical point of this branching process is known to correspond to the case in which the expectation of descendants equals $1$, see for instance \cite{grimmett2} 5.4 Theorem (5). $\square$

\begin{remark}
\textnormal{ Analyzing the equation (\ref{eq:criticaltree}), one concludes that the probability that a given site $x$ is vacant at the critical level $u_*$ is
\begin{equation}
\mathbb{P}[x \in \mathcal{V}^{u_*}] = \text{exp}\{-u_* \cdot \text{cap}(\{x\})\} = \text{exp} \left\{\frac{-d \, \text{log}(d-1)}{d-2} \right\} = \frac{1}{d}(1 + o(1)),
\end{equation}
and $\text{cap}(\{x\}) = (d-2)/(d-1)$. This gives rise to a natural question. Does the interlacement model in $\mathbb{Z}^d$ present a mean field behavior in high dimensions? Or more precisely, does the probability that the origin is vacant at the critical value of $\mathbb{Z}^d$ asymptotically behave like $1/2d$ as $d$ goes to infinite? This result is known to hold true for the case of Bernoulli independent percolation, see for instance \cite{kesten}, \cite{alon}.}
\end{remark}

\end{document}